\documentclass[12pt,times,review,sort&compress]{elsarticle}
\usepackage{epstopdf}
\usepackage[caption=false]{subfig}
\usepackage[hidelinks]{hyperref}
\usepackage[active]{srcltx}
\usepackage{amsthm}
\usepackage{amssymb}
\usepackage{stackrel}
\usepackage{xcolor}
\usepackage{amsmath,amscd}

\newtheorem{dl}{Theorem}[section]
\newtheorem{tl}[dl]{Corollary}
\newtheorem{yl}[dl]{Lemma}
\newtheorem{dy}[dl]{Definition}

\newtheorem{prop}[dl]{Proposition}

\newtheorem{remark}[dl]{Remark}

\numberwithin{equation}{section}
\newcommand{\Ome}{\mathbf{\Omega}}

\newcommand{\be}{\begin{equation}}
\newcommand{\ee}{\end{equation}}
\newcommand{\ba}{\begin{array}}
\newcommand{\ea}{\end{array}}
\newcommand{\bmn}{\begin{eqnarray}}
\newcommand{\emn}{\end{eqnarray}}
\newcommand{\bnm}{\begin{eqnarray*}}
\newcommand{\enm}{\end{eqnarray*}}
\newcommand{\bln}{\begin{subequations}}
\newcommand{\eln}{\end{subequations}}

\newproof{pot333}{Proof of Theorem \ref{firsteq}}
\newproof{pot444}{Proof of Theorem \ref{secondeq-second}}
\newproof{pot777}{Proof of Theorem \ref{secondeq-third}}

\newproof{pot555}{Applications of Theorem \ref{firsteq}}
\newproof{pot666}{Applications of Theorem \ref{secondeq-second}}
\newproof{pot888}{Applications of Theorem \ref{secondeq-third}}
\newcommand{\poq}[2]{(#1;q)_{#2}}

\def\qed{\hfill \rule{4pt}{7pt}}
\def\pf{\noindent {\it Proof.} }

\begin{document}

\title{Three new $q$-Abel  transformations and their applications}
\author{Jianan Xu\fnref{fn1}}
\fntext[fn1]{E-mail address: 20224007007@stu.suda.edu.cn}
\address[P.R.China]{Department of Mathematics, Soochow University, SuZhou 215006, P.R.China}
\author{Xinrong Ma\fnref{fn3}}
\fntext[fn3]{Corresponding author. E-mail address: xrma@suda.edu.cn.}
\address[P.R.China]{Department of Mathematics, Soochow University, SuZhou 215006, P.R.China}
\markboth{J.N. Xu and  X. R. Ma}{On two $q$-Abel series transformations and their applications}
\begin{abstract}
In the present paper, we establish three special  $q$-Abel transformation formulae of $q$-series via the use of  Abel's lemma on summation by parts. As direct applications, we set up the corresponding $q$-contiguous relations for three kinds of truncated  $q$-series.  Several new transformations are consequently established.\end{abstract}
\maketitle
\vspace{20pt}
\parskip 7pt
\section{Introduction}
It is well known  that Abel's lemma on summation by parts is one of the most  important tools in  Analysis and Number Theory, which is analogous to the integration by parts in Calculus. Abel's lemma can often be stated as follows:
\begin{yl}[Abel's lemma on summation by parts: {\rm \cite[p. 313]{knopp}}]\label{abelemma} For integer $n\geq 0$ and sequences
$\{A_n\}_{n\geq 0}$ and $\{B_n\}_{n\geq 0}$, it always holds
\begin{align}\sum_{k=0}^{n-1} B_k\Delta A_k
=A_{n-1}B_{n}-A_{-1}B_0+\sum_{k=0}^{n-1}A_k\nabla B_k.\label{abelparts}
\end{align}
Hereafter, we define the backward and forward difference operators $\Delta$ and $\nabla$ acting on arbitrary sequence $\{X_n\}_{n\geq 0}$, respectively, by
\begin{align}
\Delta X_k:=X_k-X_{k-1}\quad and \quad
\nabla X_k:=X_k-X_{k+1}.
\end{align}
\end{yl}
To  make this paper self-contained, we would like to reiterate some background from \cite{xuxrma}. First of all,  it should be mentioned that it is Abel's lemma on summation by parts with which W. C. Chu and his cooperators via a series of papers,
 say \cite{chen0,chen,chu06, chu2007, Chu1, chu2008,chu2009},
 have systematically investigated  and showed many transformation and summation formulae from the theory of various  (general, basic, and theta) hypergeometric series. Among these proofs,  perhaps  most noteworthy is that by Abel's lemma,  Chu \cite{chu06} found an elementary proof for the well-known Bailey ${}_6\psi_6$ summation formula.  To the best of our knowledge,    Chu et al.'s  elementary proofs, to great extent, depend on  two basic algebraic identities.

\begin{yl}\label{lemm12}For any complex parameters $a,b,c,x,y,z,w: xyzw\neq 0$, we define
\begin{align}\theta(x;q):=(x;q)_\infty(q/x;q)_\infty,~~\mathcal{D}(x):=1-x.\label{charafunchi}
\end{align}
Then the following two identities are true.
\begin{enumerate} [(i)]
\item (Four-term algebraic identity)
 \begin{align}
 \mathcal{D}\big(cx,\frac{x}{c},bz,\frac{z}{b}\big)-\mathcal{D}\big(bx,\frac{x}{b},cz,\frac{z}{c}\big)=\frac{z}{c}\mathcal{D}\big(bc,\frac{c}{b},xz,\frac{x}{z}\big).\label{kkklll-111}
 \end{align}
\item (Weierstrass' theta identity: \cite[Ex. 20.5.6]{10})
 \begin{align}
 \theta\big(cx,\frac{x}{c},bz,\frac{z}{b};q\big)-\theta\big(bx,\frac{x}{b},cz,\frac{z}{c};q\big)=\frac{z}{c}\theta\big(bc,\frac{c}{b},xz,\frac{x}{z};q\big).\label{trivalweierstrass}
 \end{align}
\end{enumerate}
In the above, the multivariate notation
$\mathcal{D}(x_1,x_2,\ldots,x_m):\:=\:\mathcal{D}(x_1)\mathcal{D}(x_2)
\ldots\mathcal{D}(x_m)$  and
\(\theta(x_1,x_2,\ldots,x_m):\:=\:\theta(x_1)\theta(x_2)\ldots\theta(x_m),\) $m\geq 1$ being  integer.
\end{yl}
It is worthy pointing out that   \eqref{kkklll-111}  is equivalent to Weierstrass' theta identity \eqref{trivalweierstrass}. As for this equivalency, we refer the reader to \cite{wangjin} for  a full proof. In addition, we would like to refer the reader to \cite{koornwinder} for a historical  anecdote about the origins of Weierstrass' theta identity.

More or less, the idea of applying Abel's lemma to $q$-series may  trace back to Gasper and Rahman's works  \cite{Gasper89}{} and \cite{Gasper90},  wherein they  set forth  a lot of multibasic  (quadratic, cubic, and quartic) transformations for $q$-series via the difference method which is a specific  Abel's lemma.  Particularly, Gasper believed that their method could not work in general setting, the reason is that for arbitrary sequence
\begin{align}
s_n=\frac{(a ; p)_n(b ; q)_n(c ; P)_n(d ; Q)_n}{\big(e ; p^{\prime}\big)_n\big(f ; q^{\prime}\big)_n\big(g; P^{\prime}\big)_n\big(h ; Q^{\prime}\big)_n},
\end{align}
 the corresponding difference $\Delta(s_n)=s_n-s_{n-1}$ does not take the factorial  form. The reader might consult \cite[Eq. (2.8)]{Gasper89}
for further  details. Unlike Gasper's intuition,  there do exist some  general sequences $s_n$ such that $\Delta(s_n)$ can be decomposed into product of factorial  factors, provided that  $s_n$  is subject to \eqref{kkklll-111}.    Having this fact in hand,  we have made in \cite{xuxrma}  a systematic study of Gasper's multibasic transformations
and finally established the following rather general transformation, which is now coined  {\sl  the general $q$-Abel transformation}
in order to emphasis the central role  played by Abel's lemma.
\begin{yl}[{\rm The general $q$-Abel transformation: \cite[Thm.2.1]{xuxrma}}]\label{type-i-i}
For any integer $n\geq 0$ and any complex vectors $\bar{a}=(a_1,a_2,a_3,a_4),$
$\bar{p}=(p_1,p_2,p_3,p_4)$, $\bar{x}=(x_1,x_2,x_3,x_4), \bar{q}=(q_1,q_2,q_3,q_4)$, where $a_ix_ip_iq_i\neq 0,1\leq i\leq 4$, there holds
\begin{align}
&\sum_{k=0}^{n-1}\Gamma_{k-1}[\bar{a};\bar{p}]\frac{KL^{k-1}-1}{KL^{k-1} }~ \prod_{i=1}^4\frac{(a_i^2;p_i)_{k-1}}{(K/a_i^2;L/p_i)_k
}\prod_{i=1}^4
\frac{(x_i^2;q_i)_k}
{(K_0/x_i^2;L_0/q_i)_k}\nonumber\\
  &=
\prod_{i=1}^4\frac{(a_i^2;p_i)_{n-1}}
{(K/a_i^2;L/p_i)_{n-1}}\prod_{i=1}^4\frac{(x_i^2;q_i)_n}
{(K_0/x_i^2;L_0/q_i)_{n}}-
\prod_{i=1}^4\frac
{L/p_i-K/a_i^2}{p_i-a_i^2}\label{pppppp}\\
&+\sum_{k=0}^{n-1}\Gamma_{k}[\bar{x};\bar{q}]\frac{1-K_0L_0^{k}}{K_0L_0^{k} }
~ \prod_{i=1}^4\frac{(x_i^2;q_i)_k}
{(K_0/x_i^2;L_0/q_i)_{k+1}}
\prod_{i=1}^4\frac{(a_i^2;p_i)_k}
{(K/a_i^2;L/p_i)_k}\nonumber
.
\end{align}
Hereafter, we define
\begin{align}K:=a_1a_2a_3a_4,\quad& L:=(p_1p_2p_3p_4)^{1/2},\label{kkklll-1}\\
K_0:=x_1x_2x_3x_4,\quad&L_0:=(q_1q_2q_3q_4)^{1/2},\label{kkklll-2}\end{align}
and for integer $k\geq -1$, the function
\begin{align}
\Gamma_k[\bar{x};\bar{q}]&:= \big(x_1x_2 (q_1q_2)^{k/2}-x_3x_4 (q_3q_4)^{k/2}\big) \label{chachacha}\\
&\times\big(x_1x_3 (q_1q_3)^{k/2}-x_2x_4
(q_2q_4)^{k/2}\big)\big(x_1x_4 (q_1q_4)^{k/2}-x_2x_3 (q_2q_3)^{k/2}\big).\nonumber
\end{align}
\end{yl}

Evidently,  the problem of  studying applications of \eqref{pppppp} becomes necessary and attractive. In this regard,   we refer the reader to our previous paper \cite{xuxrma} for some relevant preliminary results.  In the present paper, as further development of \cite{xuxrma}, we will derive  three  special $q$-Abel transformations, i.e., Theorems \ref{firsteq}, \ref{secondeq-second}, and \ref{secondeq-third} below, from Lemma  \ref{type-i-i} with $\bar{a}$ and $\bar{x}$ subject to three different substitutions.
In the sequel, we will focus on applications of these special $q$-Abel transformations to the theory of $q$-series.

Our paper is organized as follows. In the next   section, we will  show how to derive  Theorems \ref{firsteq}, \ref{secondeq-second}, and \ref{secondeq-third} from Lemma \ref{type-i-i}.  Section 3 is devoted to some concrete   transformations which are deducible from  Theorems \ref{firsteq}, \ref{secondeq-second}, and \ref{secondeq-third}. Among include  some new $q$-contiguous relations  and transformations for $q$-series which seem  to have not been known in the literature. We believe that these $q$-contiguous relations reflect the underlying phenomenon of Abel's lemma. In this sense,  an open problem is proposed in Section 4.

Some explanations on notation are necessary. Throughout this paper, we will follow the notation and terminology in the book \cite{10} by Gasper and Rahman. The $q$-shifted factorial of complex variable $z$ with the base $q: |q|<1$, is given by
\begin{align*}
(a;q)_{\infty}:=\prod_{k=0}^{\infty}(1-aq^k), ~~(a;q)_{n}:=\frac{(a;q)_{\infty}}{(aq^n;q)_{\infty}}~~(n\in \mathbf{Z}),
\end{align*}
where $\mathbf{Z}$ denotes the set of integers.
For  any complex numbers $a_1,a_2,\ldots,a_m$, the multivariate notation
\begin{align*}
(a_1,a_2,\ldots,a_m;q)_n:=\prod_{k=1}^{m}(a_k;q)_n
\end{align*}
and the basic hypergeometric series with the base $q$ and the variable $z$ is defined to be
\begin{align}
_{r}\phi_{r-1}\left[\begin{matrix}a_1,a_2,\ldots,a_r\\b_1,b_2,\ldots,b_{r-1}\end{matrix};q,z\right]
:=\sum_{n=0}^{+\infty}\frac{(a_1,a_2,\ldots,a_r;q)_n}{(b_1,b_2,\ldots,b_{r-1};q)_n}\frac{z^{n}}{\poq{q}{n}}.\label{phiseries-1}
\end{align}
\section{The main theorems and proofs}
This section  is  devoted to  the main theorems of the present paper  and their  full proofs.

\begin{dl}[The first special $q$-Abel transformation]\label{firsteq}With the same notation as Lemma \ref{type-i-i}. Then, for any integer $n\geq 0$, we have
 \begin{align}
\sum_{k=0}^{n-1}\frac{E_{k-1}[\bar{a};\bar{p}]}{(p_3p_4)^{(k-1)/2}}
\frac{(a_1^2/p_1;p_1)_{k}(a_2^2/p_2;p_2)_{k}}
{(a_1a_2p_3X/L;L/p_3)_{k+1}(a_1a_2p_4/LX;L/p_4)_{k+1}
}\nonumber\\
\qquad\qquad\qquad\qquad\times\frac{(x_1^2;q_1)_k(x_2^2;q_2)_k}{(x_1x_2Y;L_0/q_3)_{k}(x_1x_2/Y;L_0/q_4)_{k}}\nonumber\\
 =1-\frac{(a_1^2/p_1;p_1)_{n}(a_2^2/p_2;p_2)_{n}(x_1^2;q_1)_n(x_2^2;q_2)_n}{
(a_1a_2p_3X/L;L/p_3)_{n}(a_1a_2p_4/LX;L/p_4)_{n}
(x_1x_2Y;L_0/q_3)_n(x_1x_2/Y;L_0/q_4)_n}\label{oooooo-1}\\
-\sum_{k=0}^{n-1}\frac{E_k[\bar{x};\bar{q}]}{(q_3q_4)^{k/2}} \frac{(a_1^2/p_1;p_1)_{k+1}(a_2^2/p_2;p_2)_{k+1}}
{(a_1a_2p_3X/L;L/p_3)_{k+1}(a_1a_2p_4/LX;L/p_4)_{k+1}
}\nonumber\\
\qquad\qquad\qquad\qquad\times\frac{(x_1^2;q_1)_k(x_2^2;q_2)_k}{
(x_1x_2Y;L_0/q_3)_{k+1}(x_1x_2/Y;L_0/q_4)_{k+1}},
\nonumber
\end{align}
where the notation $X:=a_4/a_3, Y:=x_4/x_3$, and for any integer $k\geq -1$, we define
\begin{align}
E_k[\bar{x};\bar{q}]&:=\big(x_1(q_1q_3)^{k/2}-x_2x_4
(q_2q_4)^{k/2}/x_3\big)~\big(x_1(q_1q_4)^{k/2}-x_2x_3
(q_2q_3)^{k/2}/x_4\big).\label{YYY-123}
\end{align}
\end{dl}
\pf To show this theorem,  we need only to set
$$(a_3,a_4)\to (a_3s,a_4s)\quad\mbox{and}\quad(x_3,x_4)\to(x_3t,x_4t)$$
in \eqref{pppppp}. Under such a replacement of  parameters, we find that $L$ and $L_0$ remain unchanged but $K\to Ks^2$ and $K_0\to K_0t^2$. Consequently, the functions  \[\Gamma_{k-1}[\bar{a};\bar{p}]=\Gamma_{k-1}[\bar{a};\bar{p}]_{s}s^2\quad\mbox{and}\quad \Gamma_k[\bar{x};\bar{q}]=\Gamma_k[\bar{x};\bar{q}]_tt^2,\]
where
\begin{align*}
\Gamma_{k-1}[\bar{a};\bar{p}]_{s}&:= \big(a_1a_2 (p_1p_2)^{(k-1)/2}-a_3a_4 s^2 (p_3p_4)^{(k-1)/2}\big)\mbox{Temp}_{k-1}[\bar{a};\bar{p}]
\end{align*}
and
\begin{align*}
\Gamma_k[\bar{x};\bar{q}]_{t}&:= \big(x_1x_2 (q_1q_2)^{k/2}-x_3x_{4}t^2 (q_3q_4)^{k/2}\big)\mbox{Temp}_k[\bar{x};\bar{q}]\nonumber
\end{align*}
by defining
 \begin{align*}
\mbox{Temp}_{k}[\bar{x};\bar{q}]&:=\big(x_1x_{3}(q_1q_3)^{k/2}-x_2
x_{4}(q_2q_4)^{k/2}\big)\big(x_1x_{4}(q_1q_4)^{k/2}-x_2 x_{3}(q_2q_3)^{k/2}\big).
\end{align*}
With the help of these results, we may specialize \eqref{pppppp} to
\begin{align}
S_1=S_2+S_3,\label{S123}
\end{align}
where $S_i~(i=1,2,3)$ are define, respectively, by
\begin{align*}
S_1:=\sum_{k=0}^{n-1}\Gamma_{k-1}[\bar{a};\bar{p}]_{s}{s^2} ~
\frac{Ks^2L^{k-1}-1}{K{s^2}L^{k-1} }\\
\qquad\times  \frac{(a_1^2;p_1)_{k-1}(a_2^2;p_2)_{k-1}(a_3^2s^2;p_3)_{k-1}
(a_4^2s^2;p_4)_{k-1}}{(Ks^2/a_1^{2};L/p_1)_k(Ks^2/a_2^{2};L/p_2)_k(K/a_3^{2};L/p_3)_k(K/a_{4}^2;L/p_4)_k
}\\
\qquad\times
\frac{(x_1^2;q_1)_k(x_2^2;q_2)_k(x_3^2t^2;q_3)_k(x_4^2t^2;q_4)_k}
{(K_{0}t^2/x_1^{2};L_0/q_1)_{k}(K_{0}t^2/x_2^{2};L_0/q_2)_{k}
(K_{0}/x_{3}^2;L_0/q_3)_{k}(K_{0}/x_{4}^{2};L_0/q_4)_{k}};
\end{align*}
\begin{align*}
 S_2 :=-
\frac
{L/p_1-Ks^2/a_1^{2}}{p_1-a_1^2}\frac
{L/p_2-Ks^2/a_2^{2}}{p_2-a_2^2}\frac
{L/p_3-K/a_{3}^{2}}{p_3-a_{3}^{2}s^2}\frac
{L/p_4-K/a_{4}^{2}}{p_4-a_{4}^{2}s^2}\\
+\frac{(a_1^2;p_1)_{n-1}(a_2^2;p_2)_{n-1}(a_3^2s^2;p_3)_{n-1}
(a_4^2s^2;p_4)_{n-1}}{(Ks^2/a_1^{2};L/p_1)_{n-1}(Ks^2/a_2^{2};L/p_2)_{n-1}
(K/a_{3}^{2};L/p_3)_{n-1}(K/a_4^2;L/p_4)_{n-1}
}\\
\times\frac{(x_1^2;q_1)_n(x_2^2;q_2)_n(x_3^2t^2;q_3)_n(x_4^2t^2;q_4)_n}
{(K_{0}t^2/x_1^{2};L_0/q_1)_{n}(K_{0}t^2/x_2^{2};L_0/q_2)_{n}
(K_{0}/x_{3}^2;L_0/q_3)_{n}(K_{0}/x_{4}^{2};L_0/q_4)_{n}};
\end{align*}
\begin{align*}
S_3:=\sum_{k=0}^{n-1}\Gamma_k[\bar{x};\bar{q}]_{t}{t^2}~ \frac{1-K_{0}t^2L_0^{k}}
{K_{0}{t^2}L_0^{k} }\\
\times \frac{(a_1^2;p_1)_{k}(a_2^2;p_2)_{k}(a_3^2s^2;p_3)_{k}
(a_4^2s^2;p_4)_{k}}{(Ks^2/a_1^{2};L/p_1)_k(Ks^2/a_2^{2};L/p_2)_k
(K/a_{3}^{2};L/p_3)_kK/a_{4}^{2};L/p_4)_k
}\\
\times  \frac{(x_1^2;q_1)_k(x_2^2;q_2)_k(x_3^2t^2;q_3)_k(x_4^2t^2;q_4)_k}
{(K_{0}t^2/x_1^{2};L_0/q_1)_{k+1}(K_{0}t^2/x_2^{2};L_0/q_2)_{k+1}
(K_{0}/x_{3}^2;L_0/q_3)_{k+1}(K_{0}/x_{4}^{2};L_0/q_4)_{k+1}}\nonumber.
\end{align*}
Next, by canceling  $s$ and $t$ in the denominators and taking $s\to 0, t\to 0$ simultaneously, we simplify \eqref{S123}  to
\begin{align*}
&-\sum_{k=0}^{n-1}
\frac{\Gamma_{k-1}[\bar{a};\bar{p}]_{s=0}}{KL^{k-1} }\frac{(a_1^2;p_1)_{k-1}(a_2^2;p_2)_{k-1}(x_1^2;q_1)_k(x_2^2;q_2)_k}{(K/a_3^{2};L/p_3)_k(K/a_{4}^2;L/p_4)_k
(K_{0}/x_{3}^2;L_0/q_3)_{k}(K_{0}/x_{4}^{2};L_0/q_4)_{k}}\nonumber\\
  &=-\frac
{L/p_3-K/a_{3}^{2}}{p_1-a_1^2}\frac
{L/p_4-K/a_{4}^{2}}{p_2-a_2^2}\\
&+ \frac{(a_1^2;p_1)_{n-1}(a_2^2;p_2)_{n-1}(x_1^2;q_1)_n(x_2^2;q_2)_n}{
(K/a_{3}^{2};L/p_3)_{n-1}(K/a_4^2;L/p_4)_{n-1}
(K_{0}/x_{3}^2;L_0/q_3)_n(K_{0}/x_{4}^{2};L_0/q_4)_n}\\
&+\sum_{k=0}^{n-1}\frac{\Gamma_k[\bar{x};\bar{q}]_{t=0}}
{K_{0}L_0^{k} }\frac{(a_1^2;p_1)_{k}(a_2^2;p_2)_{k}(x_1^2;q_1)_k(x_2^2;q_2)_k}{(K/a_{3}^{2};L/p_3)_k(K/a_{4}^{2};L/p_4)_k(K_{0}/x_{3}^2;L_0/q_3)_{k+1}(K_{0}/x_{4}^{2};L_0/q_4)_{k+1}}
\nonumber.
\end{align*}
A substitution of $\Gamma_{k-1}[\bar{a};\bar{p}]_{s=0}$ and $\Gamma_k[\bar{x};\bar{q}]_{t=0}$ yields
\begin{align*}
&-\frac{1}{a_3a_4}\sum_{k=0}^{n-1}
\frac{\mbox{Temp}_{k-1}[\bar{a};\bar{p}]}{(p_3p_4)^{(k-1)/2}}~ \frac{(a_1^2;p_1)_{k-1}(a_2^2;p_2)_{k-1}(x_1^2;q_1)_k(x_2^2;q_2)_k}{(K/a_3^{2};L/p_3)_k(K/a_{4}^2;L/p_4)_k
(K_{0}/x_{3}^2;L_0/q_3)_{k}(K_{0}/x_{4}^{2};L_0/q_4)_{k}}\nonumber\\
  &=-\frac{1-Kp_3/a_{3}^{2}L}{1-a_1^2/p_1}\frac
{1-Kp_4/a_{4}^{2}L}{1-a_2^2/p_2}\\
&\qquad\qquad+ \frac{(a_1^2;p_1)_{n-1}(a_2^2;p_2)_{n-1}(x_1^2;q_1)_n(x_2^2;q_2)_n}{
(K/a_{3}^{2};L/p_3)_{n-1}(K/a_4^2;L/p_4)_{n-1}
(K_{0}/x_{3}^2;L_0/q_3)_n(K_{0}/x_{4}^{2};L_0/q_4)_n}\\
&+\frac{1}{x_3x_4}\sum_{k=0}^{n-1}\frac{\mbox{Temp}_k[\bar{x};\bar{q}]}
{(q_3q_4)^{k/2} }~\frac{(a_1^2;p_1)_{k}(a_2^2;p_2)_{k}(x_1^2;q_1)_k(x_2^2;q_2)_k}{(K/a_{3}^{2};L/p_3)_k(K/a_{4}^{2};L/p_4)_k
(K_{0}/x_{3}^2;L_0/q_3)_{k+1}(K_{0}/x_{4}^{2};L_0/q_4)_{k+1}}\nonumber
.
\end{align*}
As  a last step, we reformulate this identity in terms of $X, Y, E_k[\bar{x};\bar{q}]$ and multiple both sides with
\[-\frac{1-a_1^2/p_1}{1-Kp_3/a_{3}^{2}L}\frac
{1-a_2^2/p_2}{1-Kp_4/a_{4}^{2}L}.\]
Hence the theorem is proved.
\qed

\begin{dl}[The second special $q$-Abel transformation]\label{secondeq-second}With the same notation as Lemma \ref{type-i-i}. Then, for any integer $n\geq 0$, we have
\begin{align} &\sum_{k=0}^{n-1}
\frac{E_{k-1}[\bar{a};\bar{p}]}{(p_3p_4)^{(k-1)/2}}\frac{(a_1^2/p_1;p_1)_{k}(a_2^2/p_2;p_2)_{k}(x_1^2;q_1)_k}{(Kp_3/La_3^2;L/p_3)_{k+1}(Kp_4/La_4^2;L/p_4)_{k+1}
(K_0/x_2^2;L_0/q_2)_{k}}\nonumber\\
  &=1- \frac{(a_1^2/p_1;p_1)_{n}(a_2^2/p_2;p_2)_{n}(x_1^2;q_1)_n}{
(Kp_3/La_3^2;L/p_3)_{n}(Kp_4/La_4^2;L/p_4)_{n}
(K_0/x_2^2;L_0/q_2)_n} \\
&-\frac{x_1}{x_2}\sum_{k=0}^{n-1}D_k[\bar{x};\bar{q}]\bigg(\frac{q_1}
{q_2}\bigg)^{k/2}\frac{(a_1^2/p_1;p_1)_{k+1}(a_2^2/p_2;p_2)_{k+1}(x_1^2;q_1)_k}{(Kp_3/La_3^2;L/p_3)_{k+1}(Kp_4/La_4^2;L/p_4)_{k+1}(K_0/x_2^2;L_0/q_2)_{k+1}},
\nonumber
\end{align}
where $E_k$ is the same as \eqref{YYY-123} and
\begin{align}
D_k[\bar{x};\bar{q}]:=x_1x_2(q_1q_2)^{k/2}-x_3x_4(q_3q_4)^{k/2}.\label{YYY-123123}
\end{align}
\end{dl}
\pf To prove Theorem \ref{secondeq-second}, we make the replacement of parameters
$$
\big(a_3, a_4\big) \rightarrow\big(a_3 s, a_4 s\big)~~ \text { and }~~(x_2, x_3, x_4) \rightarrow(x_2 t^2, x_3 t, x_4 t)
$$
in \eqref{pppppp}. As have made in the proof of Theorem \ref{firsteq}, by simplifying and taking the limit as  $s, t$ tend to zero, we
have the desired  transformation. The other detail is similar to above and thus left to the reader.
\qed

\begin{dl}[The third special $q$-Abel transformation]\label{secondeq-third}With the same notation as Lemma \ref{type-i-i} and Theorem \ref{secondeq-second}. Then we have
\begin{align} &\frac{a_1}{a_2}\sum_{k=0}^{n-1}D_{k-1}[\bar{a};\bar{p}]
\bigg(\frac{p_1}{p_2}\bigg)^{(k-1)/2}
\frac{(a_1^2/p_1;p_1)_{k}(x_1^2;q_1)_k}{
(Kp_2/La_2^2;L/p_2)_{k+1}
(K_0/x_2^2;L_0/q_2)_{k}}\nonumber\\
  &=1- \frac{(a_1^2/p_1;p_1)_{n}(x_1^2;q_1)_n}{
(Kp_2/La_2^2;L/p_2)_{n}
(K_0/x_2^2;L_0/q_2)_n} \label{eq116}\\
&-\frac{x_1}{x_2}\sum_{k=0}^{n-1}D_k[\bar{x};\bar{q}]
\bigg(\frac{q_1}
{q_2}\bigg)^{k/2}\frac{(a_1^2/p_1;p_1)_{k+1}
(x_1^2;q_1)_k}{(Kp_2/La_2^2;L/p_2)_{k+1}(K_0/x_2^2;L_0/q_2)_{k+1}},
\nonumber
\end{align}
\end{dl}
\pf To obtain this theorem,  we apply the change of parameters
$$
(a_2,a_3, a_4) \rightarrow(a_2 s^2, a_3 s, a_4 s)~~ \text { and }~~(x_2, x_3, x_4) \rightarrow(x_2 t^2, x_3 t, x_4 t)
$$
to \eqref{pppppp} and carry out the same procedure as above.  We leave the detail to the interested reader.
\qed
\begin{remark}We think that  Theorems   \ref{firsteq} -- \ref{secondeq-third} are completely different from each other, the distinction consist in   Theorem    \ref{secondeq-third}  can not be deduced from Theorem  \ref{secondeq-second} just only  letting $a_2=0$ while Theorem    \ref{secondeq-second}  is not the special case $x_2=0$ of  Theorem  \ref{firsteq}. Later as we will see, they lead us to different $q$-series transformations.
\end{remark}
\section{Applications}Throughout the sequel, we will establish  a few  special cases of Theorems   \ref{firsteq} -- \ref{secondeq-third} and then investigate their applications to $q$-series.

\subsection{Applications of Theorem \ref{firsteq}}
To begin, we  need to establish a specific but useful case   of Theorem \ref{firsteq}.

\begin{dl}\label{firsteqcorl}With the same notation as Theorem \ref{firsteq}. Let $ p_1=(q_1 q_2 q_4/q_3)^{1/2}$. Then it holds for any integer $n\geq 0$
\begin{align}
&\sum_{k=0}^{n-1}
\big((p_2 p_4/p_3)^{(k-1)/2}a_2X-p_1^{(k-1)/2}a_1\big)
\big((p_2 p_3/p_4)^{(k-1)/2}a_2/X-p_1^{(k-1)/2}a_1\big)\nonumber\\
&\quad\times~ \frac{(a_2^2/p_2;p_2)_{k}(x_1^2;q_1)_k(x_2^2;q_2)_k}{(a_1a_2p_3X/L;L/p_3)_{k+1}
(a_1a_2p_4/XL;L/p_4)_{k+1}
(x_1^2x_2^2p_1/a_1^2;L_0/q_4)_{k}}\nonumber\\
  &=1-\frac{(a_2^2/p_2;p_2)_{n}(x_1^2;q_1)_n(x_2^2;q_2)_n}{
(a_1a_2p_3X/L;L/p_3)_{n}(a_1a_2p_4/XL;L/p_4)_{n}
(x_1^2x_2^2p_1/a_1^2;L_0/q_4)_n}\label{threethreeeq}\\
&-\frac{1-a_2^2/p_2}{1-x_1^2x_2^2p_1/a_1^2}
\sum_{k=0}^{n-1}
\big((q_1 q_4/q_3)^{k/2} -q_2^{k/2}p_1x_2^2/a_1^2\big)\big((q_1 q_3/ q_4)^{k/2}x_1^2-q_2^{k/2}a_1^2/p_1\big)\nonumber\\
&\quad\times~\frac{(a_2^2;p_2)_{k}(x_1^2;q_1)_k(x_2^2;q_2)_k}{(a_1a_2p_3X/L;L/p_3)_{k+
1}(a_1a_2p_4/XL;L/p_4)_{k+1}
(x_1^2x_2^2p_1L_0/q_4a_1^2;L_0/q_4)_{k}}.
\nonumber
\end{align}
\end{dl}
\pf  To show \eqref{threethreeeq}, we need to assume further that, besides $p_1=L_0/q_3=(q_1 q_2 q_4/q_3)^{1/2}$,
\begin{align}a_1^2/p_1=x_1x_2Y,
 \quad\mbox{i.e.,}\quad Y=a_1^2/p_1x_1x_2. \label{conditionY}
\end{align}
As such,  we specialize Theorem \ref{firsteq}  to
\begin{align*}&\sum_{k=0}^{n-1}
\frac{E_{k-1}[\bar{a};\bar{p}]}{(p_3p_4)^{(k-1)/2}}~ \frac{(a_2^2/p_2;p_2)_{k}(x_1^2;q_1)_k(x_2^2;q_2)_k}{(a_1a_2p_3X/L;L/p_3)_{k+1}
(a_1a_2p_4/XL;L/p_4)_{k+1}
(x_1^2x_2^2p_1/a_1^2;L_0/q_4)_{k}}\nonumber\\
  &=1-\frac{(a_2^2/p_2;p_2)_{n}(x_1^2;q_1)_n(x_2^2;q_2)_n}{
(a_1a_2p_3X/L;L/p_3)_{n}(a_1a_2p_4/XL;L/p_4)_{n}
(x_1^2x_2^2p_1/a_1^2;L_0/q_4)_n}-\frac{1-a_2^2/p_2}{1-x_1x_2/Y}\nonumber\\
&\times\sum_{k=0}^{n-1}\frac{E_k[\bar{x};\bar{q}]}
{(q_3q_4)^{k/2} }~\frac{(a_2^2;p_2)_{k}(x_1^2;q_1)_k(x_2^2;q_2)_k}{(a_1a_2p_3X/L;L/p_3)_{k+1}
(a_1a_2p_4/XL;L/p_4)_{k+1}
(x_1x_2L_0/q_4Y;L_0/q_4)_{k}}.
\nonumber
\end{align*}
Next, referring to \eqref{YYY-123} and using the  conditions that  $a_1=(x_1x_2p_1Y)^{1/2}$ and $p_1=(q_1 q_2 q_4/q_3)^{1/2}$, we may compute
\begin{align*}
E_k[\bar{a};\bar{p}]&=
(p_3p_4)^{k/2}\big(a_1p_1^{k/2}-a_2
X(p_2p_4/p_3)^{k/2}\big)~
\big(a_1p_1^{k/2}-a_2 (p_2p_3/p_4)^{k/2}/X\big);\\
E_k[\bar{x};\bar{q}]&=
(q_3q_4)^{k/2}
\big((q_1 q_4/q_3)^{k/2} -q_2^{k/2}p_1x_2^2/a_1^2\big)\big((q_1 q_3/ q_4)^{k/2}x_1^2-q_2^{k/2} a_1^2/p_1\big).
\end{align*}
Rewriting the above we get  \eqref{threethreeeq}.
\qed

Regarding  applications of Theorems   \ref{firsteq} -- \ref{secondeq-third} to $q$-series, we also need a  general solution for any recurrence relation of order one. As we will see later,  this kind of recurrence relations occurs frequently to the contiguous relation of the truncated series  related with Abel's lemma.
\begin{yl}\label{xxx-xxx} Assume the sequence $\{F_n(W)\}_{n\geq 0}$ satisfies the functional equation
\begin{align}
F_n(W)=C(W)\sigma\big(F_n(W)\big)+D_n(W),\label{analticrec-fu}
\end{align}
where the vector $W:=(a_1,a_2,\cdots,a_r)\in\mathbb{C}^r $ and the operator (namely, the replacement of parameters) $\sigma$ acting on $\mathbb{C}^r$ is defined by
\begin{align*}
\sigma(W)&:=(\sigma(a_1),\sigma(a_2),\cdots,\sigma(a_r)); \\ \sigma^0(W)&:=W;  \sigma^k(W):=\sigma\big(\sigma^{k-1}(W)).
\end{align*}
For any function $G(W)$, we define
\begin{align*}
\sigma\Big(G(W)\Big)&:=G(\sigma(W)).
\end{align*}
 Then, for $m\geq 0$, we have
\begin{align}
F_n(W)=\sigma^m(F_n(W))\prod_{k=0}^{m-1}\sigma^k(C(W))+\sum_{k=0}^{m-1}
\sigma^k(D_n(W))\prod_{i=0}^{k-1}\sigma^i(C(W)).
\label{Fn}
\end{align}
\end{yl}
\pf  Identity \eqref{Fn} can be proved by induction on $m$. We omit the detail here and leave it to the reader.
\qed

In the rest of this subsection, we turn  to study three  truncated series by virtue of  \eqref{threethreeeq}. Later as we will see, these finite sums are embed with very good $q$-contiguous relations with respect to various parameters. These $q$-contiguous relations often lead us to some $q$-series transformations.
For convenience of writing,  we record  the following concept  from \cite{xuxrma} originally proposed by ourselves.
\begin{dy}[$(R,S)$-type transformation with degree $2m$: {\rm\cite[Def. 4.2]{xuxrma}}] The transformation  \eqref{pppppp} is said to be the $(R,S)$-type with degree $2m$, provided that $$p_i=q^{2r_i},q_i=q^{2s_i}, r_i,s_i>0;~~ m=\max\{r_i,s_i\,|1\leq i\leq 4\}$$ and $R,S$ are the cardinalities of the index sets $\{r_i|1\leq i\leq 4\}$ and  $\{s_i|1\leq i\leq 4\}$, respectively.
\end{dy}
We remark that this definition   generalizes Gasper and Rahman's terminology \cite{Gasper90}  on the quadratic, cubic, quartic, and quintic transformations.

\subsubsection{$(2,3)$-Type of cubic transformation}
As planned,  we will apply  \eqref{threethreeeq} to investigate possible $q$-contiguous relation of the  finite  sum \eqref{macc-old} below which was discussed by the first author in her master thesis \cite[Chap. 3]{xujianan}:
\begin{align}
\mathbf{\Omega}_n(a,b,c):=\sum_{k=0}^{n-1}
q^{k} \frac{(b/q;q)_{k}(cq;q)_{2k}}{(aq^3,c/aq;q)_k(bcq^2;q^3)_k}.\label{macc-old}
\end{align}
Indeed, we can easily  show
\begin{dl}Let  $
\mathbf{\Omega}_n(a,b,c)$ be given by \eqref{macc-old}. Then for any integer $n\geq 0$, it holds the $q$-contiguous relation
\begin{align}
\mathbf{\Omega}_n(a,b,c)&=\frac{c(q-b)(aq;q)_2}
{(b-aq^3)(1-c)(aq-c)}\mathbf{\Omega}_n(a/q^2,bq,c/q)\label{macc}\\
  &+\frac{q(1-aq^2)  (1-bc/q)}{(b-aq^3)(1-c)}\left\{1-\frac{(b/q;q)_{n}(c;q)_{2n}}{
(aq^2,c/aq;q)_{n}(bc/q;q^3)_{n}
}\right\}\nonumber.
\end{align}
\end{dl}
\pf In order to obtain  \eqref{macc},
it suffices to take  in \eqref{threethreeeq} that
\begin{align*}
(p_1,p_2, p_3, p_4)\to (q^3,q,q^3, q)~~\mbox{and}~~(q_1,q_2, q_3, q_4)\to(q^2,q^2,q,q^3).
\end{align*}
Under such substitutions, it is easy to check $L=L_0=q^4$ and
\begin{align*}
E_k[\bar{a};\bar{p}]&=q^{3k}(a_1-a_2/X)(a_1q^{2k}-a_2
X),\\
E_k[\bar{x};\bar{q}]&=q^{3k}(x_1-x_2Yq^k)(x_1q^k-x_2/Y), ~~
Y=a_1^2/x_1x_2q^3.
\end{align*}
After
some simplification, we  specialize  \eqref{threethreeeq}  to
\begin{align}\displaystyle
&\frac{-a_2X\big(a_1-a_2/X\big)}{q(1-a_1a_2X/q)
 }\sum_{k=0}^{n-1}
\frac{q^{k}\big(1-a_1q^{2(k-1)}/a_2X\big)}{ 1-a_1a_2/Xq^3}~ \frac{(a_2^2/q;q)_{k}(x_1^2,x_2^2;q^2)_k}{(a_1a_2X;q)_k(a_1a_2/X;q^3)_k(x_1^2x_2^2q^3/a_1^2;q)_{k}}\nonumber\\
  &=1-\frac{(a_2^2/q;q)_{n}}{
(a_1a_2X/q;q)_{n}(a_1a_2/Xq^3;q^3)_{n}
}\frac{(x_1^2,x_2^2;q^2)_n}
{(x_1^2x_2^2q^3/a_1^2;q)_n}+\frac{x_1^2x_2^2q^3}{a_1^2}\frac{1-a_2^2/q}{1-x_1^2x_2^2q^3/a_1^2}\label{eq34}\\
&\quad\times~\sum_{k=0}^{n-1}\frac{q^k\big(1-a_1^2q^{k-3}/x_1^2\big)\big(1-a_1^2q^{k-3}/x_2^2\big)}{(1-a_1a_2X/q)
  (1-a_1a_2/Xq^3)}~\frac{(a_2^2;q)_{k}(x_1^2,x_2^2;q^2)_k}{(a_1a_2X;q)_k(a_1a_2/X;q^3)_k
(x_1^2x_2^2q^4/a_1^2;q)_{k}}.
\nonumber
\end{align}
Now we are ready to show \eqref{macc}.
 To do this, we need to make the parametric replacements in \eqref{eq34}, as below:
\begin{align*}
(x_1^2,x_2^2,a_1^2,a_2^2, X^2)\to(c,cq,acq^5,b,aq/bc).
\end{align*}
Under this substitution, it is easy to check
\begin{align*}
1-a_1q^{2(k-1)}/a_2X&=1-cq^{2k},\\
\big(1-a_1^2q^{k-3}/x_1^2\big)\big(1-a_1^2q^{k-3}/x_2^2\big)
&=(1-aq^{k+1})(1-aq^{k+2}).
\end{align*}
Simplifying \eqref{eq34} by these computational results,  we  get
\begin{align*}
&\frac{(b-aq^3)(1-c)}{q(1-aq^2)  (1-bc/q)}\sum_{k=0}^{n-1}
q^{k} \frac{(b/q;q)_{k}(cq;q)_{2k}}{(aq^3,c/aq;q)_k
(q^2bc;q^3)_k}\nonumber\\
  &=1-\frac{(b/q;q)_{n}(c;q)_{2n}}{
(aq^2,c/aq;q)_{n}(bc/q;q^3)_{n}
}
\\
&+\frac{c}{aq}\frac{(1-b/q)(1-aq)}
{(1-c/aq)(1-bc/q)}\sum_{k=0}^{n-1}q^k\frac{(b;q)_{k}(c;q)_{2k}}{(aq,c/a;q)_k(bcq^2;q^3)_k
}.
\nonumber
\end{align*}
Restated in terms of $\Ome_n(a,b,c)$. The conclusion follows.
\qed

\begin{remark}
As usual, we often refer  the case $n<+\infty$ of arbitrary identity like \eqref{macc} as to a contiguous relation while  the corresponding case $n=+\infty$  as to  a transformation.
\end{remark}

One of purposes for our pursuiting for  the contiguous relation \eqref{macc} is to establish useful $q$-series transformations. As is expected to be, by using Lemma \ref{xxx-xxx}
, one may find with easy
\begin{tl}[{\rm cf. \cite[Thm 3.3.1]{xujianan}}]For integers $m,n\geq 0$, it holds
\begin{align}
& \mathbf{\Omega}_n(a, b, c)-\mathbf{\Omega}_n\left(a / q^{2 m}, bq^m, c/ q^m\right) \frac{\left(1 / a q^2 ; q\right)_{2 m}(b / q ; q)_m}{(c / a q, 1 / c ; q)_m\left(b / a q^3 ; q^3\right)_m}\label{macc-1} \\
=& \frac{\left(1-a q^2\right)(1-q / b c)}{\left(1-a q^3 / b\right)(1-1 / c)}\left\{\mathcal{K}_m(a, b, c)-\mathcal{K}_m\left(aq^n, bq^n, cq^{2 n}\right) \frac{(b / q ; q)_n(c ; q)_{2 n}}{\left(a q^2, c / a q ; q\right)_n\left(b c / q ; q^3\right)_n}\right\},\nonumber
\end{align}
where
$$
\mathcal{K}_m(a, b, c)=\sum_{k=0}^{m-1} q^k \frac{(1 / a q ; q)_{2 k}(b / q ; q)_k}{(q / c, c / a q ; q)_k\left(b / a ; q^3\right)_k}.
$$
\end{tl}
\pf Observe that  \eqref{macc} can be restated as the form
\begin{align}
\mathbf{\Omega}_n(a,b,c)&=C(a,b,c)
\sigma\big(\mathbf{\Omega}_n(a,b,c)\big)+D_n(a,b,c),
\label{macc-1-11}
\end{align}
where the operator
\[\sigma(a,b,c):=(a/q^2,bq,c/q)\]
and  the coefficients
\begin{align*}C(a,b,c)&:=\frac{c(q-b)(aq;q)_2}
{(b-aq^3)(1-c)(aq-c)},\\
D_n(a,b,c)&:=\frac{q(1-aq^2)  (1-bc/q)}{(b-aq^3)(1-c)}\left\{1-\frac{(b/q;q)_{n}(c;q)_{2n}}{
(aq^2,c/aq;q)_{n}(bc/q;q^3)_{n}
}\right\}.\end{align*}
Hence, \eqref{macc-1} comes out by applying Lemma \ref{xxx-xxx} to \eqref{macc-1-11}.
\qed
\subsubsection{$(3,1)$-Type of cubic transformation}
Now we reconsider the finite sum
\be\label{Omega-3}
\mathbf{\chi}_n(a,b,c):=\sum_{k=0}^{n-1}q^{k}
\frac{(abc;q^3)_{k}(b,c;q)_k}
{(bc/q^2;q)_{2k+2}(aq;q)_{k}} .
\ee
It was studied in \cite[(1.1)]{chen0}.
By virtue of \eqref{threethreeeq}, we now discover a $q$-contiguous relation of $\mathbf{\chi}_n(a,b,c)$.
\begin{dl}[{\rm cf. \cite[p.793]{chen0}}]Let  $\mathbf{\chi}_n(a,b,c)$ be given by \eqref{Omega-3}. Then, for any integer $n\geq 0$, it holds the $q$-contiguous relation
\begin{align}
\mathbf{\chi}_n(a,b,c)&=
\frac{(1-abc)(1-aq^3/b)(1-aq^3/c)
}{(a q;q)_3}\mathbf{\chi}_n(aq^3,b,c)\label{threethreeeq-12345600}\\
&+\frac{aq^3}{bc(aq;q)_2}
\left\{1-\frac{(abc;q^3)_{n}(b,c;q)_n}{
(bc/q^2;q)_{2n}
(aq^3;q)_n}\right\}.
\nonumber
\end{align}
\end{dl}
\pf To show  \eqref{threethreeeq-12345600},
 we specificize \eqref{threethreeeq} by setting
\begin{align*}
(p_1,p_2,p_3,p_4)\to(q,q^3,q^2,q^2)~~\mbox{and}~~(q_1,q_2,q_3,q_4)\to(q,q,q,q).
\end{align*}
Under such circumstance, it is easy to check $L=q^4, L_0=q^2$ and then to reduce  \eqref{threethreeeq}  to
\begin{align}
&a_1^2\sum_{k=0}^{n-1}q^{k-1}
\big(1-a_2Xq^{k-1}/a_1\big)
\big(1-a_2q^{k-1}/a_1X\big)\nonumber\\
&\qquad\quad\times \frac{(a_2^2/q^3;q^3)_{k}(x_1^2,x_2^2;q)_k}
{(a_1a_2X/q^2,a_1a_2/Xq^2;q^2)_{k+1}
(x_1^2x_2^2q/a_1^2;q)_{k}}\nonumber\\
  &=1-\frac{(a_2^2/q^3;q^3)_{n}(x_1^2,x_2^2;q)_n}{
(a_1a_2X/q^2,a_1a_2/Xq^2;q^2)_{n}
(x_1^2x_2^2q/a_1^2;q)_n}\label{threethreeeq-123}\\
&-\big(1-x_2^2q/a_1^2\big)\big(x_1^2-a_1^2/q\big)\frac{1-a_2^2/q^3}{1-x_1^2x_2^2q/a_1^2}
\nonumber\\
&\qquad\quad\times\sum_{k=0}^{n-1}q^{k}\frac{(a_2^2;q^3)_{k}(x_1^2,x_2^2;q)_k}
{(a_1a_2X/q^2,a_1a_2/Xq^2;q^2)_{k+1}
(x_1^2x_2^2q^2/a_1^2;q)_{k}}.
\nonumber
\end{align}
  If $a_1a_2X=x_1^2x_2^2$ and $a_1a_2/X=x_1^2x_2^2q$, then
    \[a_1a_2=x_1^2x_2^2q^{1/2}~~\mbox{and}~~ X=q^{-1/2}.\]
In this case, it is easy to check
\begin{align*}
&\frac{(1-a_2Xq^{k-1}/a_1)
(1-a_2q^{k-1}/a_1X)}{
(x_1^2x_2^2q/a_1^2;q)_{k}}=\frac{(1-a_1a_2Xq^{k-1}/a_1^2)
(1-a_1a_2q^{k-1}/a_1^2X)}{
(x_1^2x_2^2q/a_1^2;q)_{k}}\\
&\quad\qquad=\frac{(1-x_1^2x_2^2q^{k-1}/a_1^2)
(1-x_1^2x_2^2q^{k}/a_1^2)}{
(x_1^2x_2^2q/a_1^2;q)_{k}}=\frac{(1-x_1^2x_2^2/a_1^2q)
(1-x_1^2x_2^2/a_1^2)}{
(x_1^2x_2^2/a_1^2q;q)_{k}}.
\end{align*}
Using these expressions, we may reformulate  both  sides of \eqref{threethreeeq-123} in the form
\begin{align}
&\frac{a_1^2\big(1-x_1^2x_2^2/a_1^2q\big)
\big(1-x_1^2x_2^2/a_1^2\big)}{q}\sum_{k=0}^{n-1}q^{k}
\frac{(x_1^4x_2^4/a_1^2q^2;q^3)_{k}(x_1^2,x_2^2;q)_k}
{(x_1x_2/q^2,x_1x_2/q;q^2)_{k+1}
(x_1^2x_2^2/a_1^2q;q)_{k}}\nonumber\\
  &=1-\frac{(x_1^4x_2^4/a_1^2q^2;q^3)_{n}(x_1^2,x_2^2;q)_n}{
(x_1^2x_2^2/q^2,x_1^2x_2^2/q;q^2)_{n}
(x_1^2x_2^2q/a_1^2;q)_n}-\big(1-x_2^2q/a_1^2\big)
\big(x_1^2-a_1^2/q\big)
\label{threethreeeq-123}\\
&\quad\times
\frac{1-x_1^4x_2^4/a_1^2q^2}{1-x_1^2x_2^2q/a_1^2}\sum_{k=0}^{n-1}q^{k}
\frac{(x_1^4x_2^4q/a_1^2;q^3)_{k}(x_1^2,x_2^2;q)_k}
{(x_1^2x_2^2/q^2,x_1^2x_2^2/q;q^2)_{k+1}
(x_1^2x_2^2q^2/a_1^2;q)_{k}}.
\nonumber
\end{align}
Evidently,  under  the  replacement of parameters
\begin{align*}
(a_1^2,x_1^2,x_2^2)\to(bc/aq^2,b,c),\end{align*}
\eqref{threethreeeq-123} assumes  the shape as \eqref{threethreeeq-12345600}.
\qed

Starting from \eqref{threethreeeq-12345600} and using Lemma \ref{xxx-xxx}, we recover without any difficulty
\begin{tl}[{\rm cf. \cite[p.796]{chen0}}]For integer $m\geq 0$, it holds
\begin{align}
&\mathbf{\chi}_\infty(a,b,c)
=\mathbf{\chi}_\infty(aq^{3m},b,c)\frac{\big(a b c, aq^3/ b, aq^3/ c;q^3)_m }{\big(aq;q)_{3m}} \nonumber\\
& +\frac{aq^3}{bc(1-aq)\left(1-aq^2 \right)} \sum_{k=0}^{m-1} q^{3 k}\frac{\big(
a b c, aq^3/b, aq^3/c;q^3\big)_k }{\big(
aq^3;q\big)_{3k}} \label{threethreeeq-12345600-new}\\
& -\frac{aq^3(b, c;q)_\infty (a b c ; q^3)_{\infty}}{bc(aq, b c/q^2;q)_{\infty}} \sum_{k=0}^{m-1}q^{3 k}(aq^3/ b ,aq^3/ c; q^3)_k .\nonumber
\end{align}
\end{tl}
\pf Observe first that  \eqref{threethreeeq-12345600} can be restated as the operator form
\begin{align}
\mathbf{\chi}_n(a,b,c)&=C(a,b,c)
\sigma\big(\mathbf{\chi}_n(a,b,c)\big)+D_n(a,b,c),
\label{threethreeeq-12345600-1}
\end{align}
where the operator
\[\sigma(a,b,c):=(aq^3,b,c)\]
and the coefficients
\begin{align*}C(a,b,c)&:=
\frac{(1-abc)(1-aq^3/b)(1-aq^3/c)
}{(a q;q)_3},\\
D_n(a,b,c)&:=\frac{aq^3}{bc(aq;q)_2}
\left\{1-\frac{(abc;q^3)_{n}(b,c;q)_n}{
(bc/q^2;q)_{2n}
(aq^3;q)_n}\right\}.\end{align*}
Hence, by Lemma \ref{xxx-xxx}, we may derive from \eqref{Fn} that
\begin{align*}
\mathbf{\chi}_n(a,b,c)&=\mathbf{\chi}_n(aq^{3m},b,c)
\prod_{i=0}^{m-1}\frac{(1-abcq^{3i})(1-aq^{3i+3}/b)(1-aq^{3i+3}/c)
}{(aq^{3i+1};q)_3}\\
&+\sum_{i=0}^{m-1}\frac{aq^{3i+3}}{bc(aq^{3i+1};q)_2}
\left\{1-\frac{(abcq^{3i};q^3)_{n}(b,c;q)_n}{
(bc/q^2;q)_{2n}
(aq^{3i+3};q)_n}\right\}\\
&\qquad\quad\times\prod_{j=0}^{i-1}\frac{(1-abcq^{3j})(1-aq^{3j+3}/b)(1-aq^{3i+3}/c)
}{(aq^{3j+1};q)_3}.
\end{align*}
Further  simplification leads us to
\begin{align*}
\mathbf{\chi}_n(a,b,c)
&=\mathbf{\chi}_n(aq^{3m},b,c)\frac{\big(a b c, aq^3/ b, aq^3/ c;q^3)_m }{\big(aq;q)_{3m}} \\
& +\frac{aq^3 }{bc(1-aq)\left(1-aq^2 \right)} \sum_{i=0}^{m-1} q^{3 i}\frac{\big(
a b c, aq^3/ b, aq^3/ c;q^3\big)_i }{\big(
aq^3;q\big)_{3i}} \\
& -\frac{aq^3(b, c;q)_n}{bc(b c/q^2;q)_{2n}} \sum_{i=0}^{m-1}q^{3i}\frac{(abc; q^3)_{n+i}(aq^3/ b,aq^3/ c; q^3)_i}{ (aq; q)_{n+3i+2}}.
\end{align*}
Finally, \eqref{threethreeeq-12345600-new} comes out by  taking the limit as $n\to \infty$.
\qed

\subsubsection{$(1,2)$-Type of quadratic transformation}
In  the same way as above, we continue to  investigate a new  finite  sum
\be
\mathbf{\Xi}_n(a,b,c,d)=\sum_{k=0}^{n-1}q^k
\frac{(a,b;q)_k(c;q^2)_k}
{(d,c/dq;q)_{k+1}(abq;q^2)_{k}
}. \label{Chieq}
\ee

\begin{dl}\label{thm3.6}Let  $
\mathbf{\Xi}_n(a,b,c,d)$ be given by \eqref{Chieq}. Then, for any integer $n\geq 0$, it holds the $q$-contiguous relation \begin{align}
\mathbf{\Xi}_n(a,b,c,d)&=
\frac{d(a;q)_2(c-abq)}{\big(c-adq\big)
\big(d-a\big)(1-abq)}
\mathbf{\Xi}_n(aq^2,b,c,d)\label{xrma-1}\\
  &+\frac{aq}{\big(c-adq\big)
\big(1-a/d\big)}\left\{1-
\frac{(a,b;q)_n(c;q^2)_n}
{(d,c/dq;q)_{n}
(abq;q^2)_n}\right\}\nonumber.
\nonumber
\end{align}
\end{dl}
\pf To show  \eqref{xrma-1},
it suffices to take  in \eqref{threethreeeq} that  $a_2^2=x_2/x_1Y$, and
\begin{align*}
(p_1,p_2,p_3,p_4)&\to(q,q,q,q)~~\mbox{and}~~(q_1,q_2,q_3,q_4)\to(q,q^2,q^2,q).
\end{align*}
In this case, it is easy to check $L=q^2,L_0=q^3$ and specify \eqref{threethreeeq} as the form
\begin{align}
&\sum_{k=0}^{n-1}
\big(a_2Xq^{(k-1)/2}-(x_1x_2Y)^{1/2}q^{k/2}\big)
\big(a_2q^{(k-1)/2}/X-(x_1x_2Y)^{1/2}q^{k/2}\big)\nonumber\\
&\qquad\quad\qquad\times  \frac{(a_2^2/q;q)_{k}(x_1^2;q)_k(x_2^2;q^2)_k}{(x_2X/q^{1/2},x_2/Xq^{1/2};q)_{k+1}
(x_1x_2/Y;q^2)_{k}}\nonumber\\
  &=1-\frac{(a_2^2/q;q)_{n}}{(x_2X/q^{1/2},x_2/Xq^{1/2};q)_{n}
}\frac{(x_1^2;q)_n(x_2^2;q^2)_n}
{(x_1x_2/Y;q^2)_n}\label{threethreeeq-eq-eq}\\
&-\frac{1-a_2^2/q}{1-x_1x_2/Y}
\sum_{k=0}^{n-1}
\big(x_1q^{k}-x_2Yq^k\big)\frac{\big(x_1-x_2q^{k}/Y\big)(a_2^2;q)_{k}
(x_1^2;q)_k(x_2^2;q^2)_k}{(x_2X/q^{1/2},x_2/Xq^{1/2};q)_{k+1}
(x_1x_2q^2/Y;q^2)_{k}}.
\nonumber
\end{align}
Observe that, when $a_2^2=x_2/x_1Y$,  it holds
\[(x_1-x_2q^{k}/Y)(a_2^2;q)_{k}=x_1(1-a_2^2)(a_2^2q;q)_{k}.\]
Therefore
\begin{align*}
&a_2^2\big(X/q^{1/2}-x_2/a_2^2\big)
\big(1/Xq^{1/2}-x_2/a_2^2\big)\sum_{k=0}^{n-1}
q^k\frac{(a_2^2/q;q)_{k}(x_1^2;q)_k(x_2^2;q^2)_k}{(x_2X/q^{1/2},x_2/Xq^{1/2};q)_{k+1}(x_1^2a_2^2;q^2)_{k}
}\nonumber\\
  &=1-\frac{(a_2^2/q;q)_{n}}{(x_2X/q^{1/2},x_2/Xq^{1/2};q)_{n}
}\frac{(x_1^2;q)_n(x_2^2;q^2)_n}
{(x_1^2a_2^2;q^2)_n}\nonumber\\
&-\frac{(a_2^2/q;q)_2(x_1^2-x_2^2/a_2^2)}{1-x_1^2a_2^2}
\sum_{k=0}^{n-1}q^{k}\frac{(a_2^2q;q)_{k}(x_1^2;q)_k(x_2^2;q^2)_k}
{(x_2X/q^{1/2},x_2/Xq^{1/2};q)_{k+1}(x_1^2a_2^2q^2;q^2)_{k}
}.
\nonumber
\end{align*}
Next, after making  the parametric replacements
\begin{align*}
(a_2^2,
x_1^2,
x_2^2,X)\to(
aq, b,c, dq^{1/2}/c^{1/2})
\end{align*}
and reformulating the resulted in terms of $\mathbf{\Xi}_n(a,b,c,d)$, we achieve \eqref{xrma-1} at once.
\qed

It is a bit surprise that by using the contiguous relation \eqref{xrma-1}, we may establish the following  new and peculiar $q$-identity.

\begin{tl}\begin{align}
&\frac{(aq;q)_{\infty}(b/d;q^2)_\infty}{(q^2,aq^2/d,abq;q^2)_\infty}
\sum_{k=0}^{\infty}q^k
\frac{(b;q)_k(adq;q^2)_k}
{(dq,aq;q)_{k}}\label{YYY-ZZZ-WWW}
\\
&=\bigg(1-\frac{1}{d}\bigg)\sum_{k=0}^{\infty}q^{2k}
\frac{(a,aq,b/d;q^2)_k}{(q^2,aq^2/d,abq;q^2)_k}
+\frac{1}{d}
\frac{(b;q)_\infty(adq;q^2)_{\infty}}{(dq;q)_\infty(abq;q^2)_\infty}
\sum_{k=0}^{\infty}q^{2k}
\frac{(b/d;q^2)_k}{(q^2,aq^2/d;q^2)_k}.\nonumber
\end{align}
\end{tl}
\pf Observe that  \eqref{xrma-1} can be restated as the form
\begin{align}
\mathbf{\Xi}_n(a,b,c,d)&=C(a,b,c,d)
\sigma\big(\mathbf{\Xi}_n(a,b,c,d)\big)+D_n(a,b,c,d),\label{xrma-1-11}
\end{align}
where the operator
\[\sigma(a,b,c,d):=(aq^2,b,c,d)\]
and the coefficients
\begin{align*}C(a,b,c,d)&:=
\frac{d(a;q)_2(c-abq)}{\big(c-adq\big)
\big(d-a\big)(1-abq)},\\
D_n(a,b,c,d)&:=\frac{aq}{\big(c-adq\big)
\big(1-a/d\big)}\left\{1-
\frac{(a,b;q)_n(c;q^2)_n}
{(d,c/dq;q)_{n}
(abq;q^2)_n}\right\}.\end{align*}
Hence, \eqref{YYY-ZZZ-WWW} comes out by applying Lemma \ref{xxx-xxx} to \eqref{xrma-1-11}
 and then taking the limit as $m,n\to \infty$, finally multiplying both sides of the  resulting identity with $(1-adq/c)(1-a/d)$ and then letting $c=adq$.
\qed

From \eqref{YYY-ZZZ-WWW} it follows at once by setting $d=1$
\begin{prop}\begin{align}
&\sum_{k=0}^{\infty}q^k
\frac{(b;q)_k(aq;q^2)_k}
{(q,aq;q)_{k}}=(-q;q)_\infty(bq;q^2)_\infty\sum_{k=0}^{\infty}q^{2k}
\frac{(b;q^2)_k}{(q^2,aq^2;q^2)_k}.\label{YYY-ZZZ-WWW-new}
\end{align}
\end{prop}
 Very interesting is that  \eqref{YYY-ZZZ-WWW-new} implies
\begin{prop}For any complex number $|t|<1$, we have
\begin{align}
\sum_{k=0}^{\infty}t^{k}
\frac{(aq;q^2)_k}
{(q,aq;q)_{k}}=
(-t;q)_\infty\sum_{k=0}^{\infty}
\frac{t^{2k}}{(q^2,aq^2;q^2)_k}.\label{finalresult}
\end{align}\end{prop}
For limitation of space, we leave  the full derivation of \eqref{finalresult} to the interested reader.

\subsection{Applications of Theorem \ref{secondeq-second}}
First of all, we need to show
\begin{dl}\label{secondeq}
 Write $\mathcal{E}_0=\big(a_1-a_2 a_4/a_3\big)\big(a_1-a_2 a_3/a_4\big)$. Then, for any integer $n\geq 0$,  we have
\begin{align}
&\mathcal{E}_0 \sum_{k=0}^{n-1}q^{k-1} \frac{\big(a_1^2 / q, a_2^2 / q ; q\big)_k\big(x_1^2 ; q^2\big)_k}{\big(K / a_3^2 q, K / a_4^2 q ; q\big)_{k+1}\big(x_1 x_3 x_4 / x_2 ; q^2\big)_k}  \nonumber\\
& =1-\frac{\big(a_1^2 / q, a_2^2 / q ; q\big)_n\big(x_1^2 ; q^2\big)_n}{\big(K / a_3^2 q, K / a_4^2 q ; q\big)_n\big(x_1 x_3 x_4 / x_2 ; q^2\big)_n} \label{eq16}\\
& - x_1(x_1-x_3x_4/x_2)\sum_{k=0}^{n-1} q^{2 k}\frac{\big(a_1^2 / q, a_2^2 / q ; q\big)_{k+1}\big(x_1^2 ; q^2\big)_k}{\big(K / a_3^2 q, K / a_4^2 q ; q\big)_{k+1}\big(x_1 x_3 x_4 / x_2 ; q^2\big)_{k+1}} .\nonumber
\end{align}
\end{dl}
\pf It  follows from Theorem \ref{secondeq-second} by setting $\big(p_i, q_i\big) \rightarrow (q, q^2)$ for $1\leq i\leq 4$.
\qed

One of the most important cases covered by Theorem \ref{secondeq} is as follows:
\begin{tl} Let $K=a_1a_2a_3a_4$ as before. Then it holds
\begin{align}
&\frac{q\big(a_1 a_3-a_2 a_4\big)\big(a_1 a_4-a_2 a_3\big)}{(a_3q-a_1a_2a_4)(a_4q-a_1a_2a_3)}~{}_{4}\phi_{3}\left[\begin{matrix}a_1^2 / q, a_2^2 / q,x_1,-x_1\\ K / a_3^2, K / a_4^2,-q \end{matrix};q,q\right]\nonumber\\
&=-\frac{\big(a_1^2 / q, a_2^2 / q,x_1,-x_1 ; q\big)_\infty}{\big(q,-q,K / a_3^2 q, K / a_4^2 q ; q\big)_\infty}+{}_{4}\phi_{3}\left[\begin{matrix}a_1^2 / q, a_2^2 / q,x_1/q,-x_1/q\\ K / a_3^2q, K / a_4^2q,-q \end{matrix};q,q^2\right].\label{eq17}
\end{align}
\end{tl}
\pf It is  the direct consequence of \eqref{eq16} by setting $x_1x_3x_4/x_2=q^2$ and taking the limit as $n\to \infty$.
\qed

Another important case contained by  \eqref{eq16} is the following transformation for $_{4}\phi_{3}$ series, which seems very different from Singh's transformation \cite[(III.21)]{10}.
\begin{tl}
\begin{align}
&{}_{4}\phi_{3}\left[\begin{matrix}a_1^2 / q, a_2^2 / q,x_1,-x_1\\  (a_1 a_2)^2/q,q,-q \end{matrix};q,q\right]\label{eq1818}\\
=&\frac{x_1^2-q^2}{1-q^2}~{}_{4}\phi_{3}\left[\begin{matrix}a_1^2 , a_2^2,x_1,-x_1\\   (a_1 a_2)^2/q,q^2,-q^2 \end{matrix};q,q^2\right]+\frac{\big(a_1^2, a_2^2,x_1,-x_1 ; q\big)_\infty}{\big(q, (a_1 a_2)^2/q,q,-q; q\big)_\infty}.\nonumber
\end{align}
\end{tl}
\pf Actually, \eqref{eq1818} is a special case of \eqref{eq17} under the assumption that $a_1a_2a_4=a_3q$. To see this clear, we reformulate \eqref{eq17} as the  equivalent form
\begin{align*}
&{}_{4}\phi_{3}\left[\begin{matrix}a_1^2 / q, a_2^2 / q,x_1,-x_1\\ K / a_3^2, K / a_4^2,-q \end{matrix};q,q\right]\nonumber\\
&\qquad=-\frac{\big(a_1^2 / q, a_2^2 / q,x_1,-x_1 ; q\big)_\infty}{\big(q,-q,K / a_3^2 q, K / a_4^2 q ; q\big)_\infty}~ \frac{(a_3q-a_1a_2a_4)(a_4q-a_1a_2a_3)}{q\big(a_1 a_3-a_2 a_4\big)\big(a_1 a_4-a_2 a_3\big)}\\
&\qquad+\frac{(a_3q-a_1a_2a_4)(a_4q-a_1a_2a_3)}{q\big(a_1 a_3-a_2 a_4\big)\big(a_1 a_4-a_2 a_3\big)}\bigg(1+\sum_{k=1}^\infty q^{2k}\frac{(a_1^2 / q, a_2^2 / q,x_1/q,-x_1/q;q)_{k}}{(q,K / a_3^2q, K / a_4^2q,-q;q)_k}\bigg).\nonumber
\end{align*}
After a bit routine simplification, it equals
\begin{align*}
&{}_{4}\phi_{3}\left[\begin{matrix}a_1^2 / q, a_2^2 / q,x_1,-x_1\\ K / a_3^2, K / a_4^2,-q \end{matrix};q,q\right]=-\frac{\big(a_1^2 / q, a_2^2 / q,x_1,-x_1 ; q\big)_\infty}{\big(q,-q,K / a_3^2, K / a_4^2 ; q\big)_\infty}~ \frac{a_3a_4q}{\big(a_1 a_3-a_2 a_4\big)\big(a_1 a_4-a_2 a_3\big)}\\
&\qquad\qquad\qquad\qquad\qquad\qquad+\frac{(a_3q-a_1a_2a_4)(a_4q-a_1a_2a_3)}{q\big(a_1 a_3-a_2 a_4\big)\big(a_1 a_4-a_2 a_3\big)}\nonumber\\
&\qquad\qquad\qquad+\frac{a_3a_4q^3(1-a_1^2/q)(1-a_2^2/q)(1-x_1^2/q^2)}{(1-q^2)\big(a_1 a_3-a_2 a_4\big)\big(a_1 a_4-a_2 a_3\big)}\sum_{k=0}^\infty q^{2k}\frac{(a_1^2, a_2^2,x_1,-x_1;q)_{k}}{(K / a_3^2, K / a_4^2,q^2,-q^2;q)_k}.\nonumber
\end{align*}
Taking the assumption  $a_1a_2a_4=a_3q$ into account, we thereby obtain
\begin{align*}
&{}_{4}\phi_{3}\left[\begin{matrix}a_1^2 / q, a_2^2 / q,x_1,-x_1\\  (a_3 / a_4)^2q,q,-q \end{matrix};q,q\right]\\
=&\frac{\big(a_1^2, a_2^2,x_1,-x_1 ; q\big)_\infty}{\big(q,(a_3 / a_4)^2q,q,-q; q\big)_\infty}+\frac{x_1^2-q^2}{1-q^2}~{}_{4}\phi_{3}\left[\begin{matrix}a_1^2 , a_2^2,x_1,-x_1\\  (a_3 / a_4)^2q,q^2,-q^2 \end{matrix};q,q^2\right].\nonumber
\end{align*}
Replacing  $a_3/a_4$ with $a_1a_2/q$, we complete  the proof of Transformation \eqref{eq1818}.\qed

\subsection{Applications of Theorem \ref{secondeq-third}}
We  conclude this section with a direct application of Theorem \ref{secondeq-third}. For this, we introduce
\begin{align}
\mathbf{\Pi}_n(a_1,x_1,y,z):=\sum_{k=0}^{n-1}
q^{k}
\frac{(a_1^2/q,x_1^2;q)_k}{
(a_1y,x_1z;q)_{k}}\label{finalsum}
\end{align}
provided that  both $a_1y$ and $x_1z$ are not of form $q^{-m}$, $m$ being nonnegative integer.
\begin{dl}\label{secondeq-third-app1}Let $\mathbf{\Pi}_n(a_1,x_1,y,z)$ be given by \eqref{finalsum}. Then there holds the $q$-contiguous relation
\begin{align}
\mathbf{\Pi}_n(a_1,x_1,y,z)&=\frac{x_1(x_1-z)(a_1^2-q)}{a_1(a_1-y)(1-x_1z)}
\mathbf{\Pi}_n(a_1q^{1/2},x_1,y/q^{1/2},zq)\nonumber\\
&+\frac{q-a_1y}{a_1(a_1-y)}
\left\{1-\frac{(a_1^2/q,x_1^2;q)_n}{(a_1y/q,x_1z;q)_n}\right\} \label{eq317317}
\end{align}
and  the transformation
\begin{align}
{}_2 \phi_1\left[\begin{array}{c}
a_1^2 / q, x_1^2 \\
a_1 y
\end{array} ; q, q\right]&+\frac{q-a_1y}{a_1(y-a_1)}{}_2 \phi_1\left[\begin{array}{c}
a_1^2/q, x_1^2/q \\
a_1 y/q
\end{array} ; q, q\right]=\frac{q}{a_1(y-a_1)}\frac{\big(a_1^2 / q,x_1^2; q\big)_{\infty}}
{(a_1 y,q; q)_{\infty}}.\label{xrma-1-xrma}
\end{align}
\end{dl}
\pf To prove \eqref{eq317317}, it is sufficient to set in \eqref{eq116} of Theorem  \ref{secondeq-third} that for $1\leq i\leq 4$, $p_i=q_i=q,$
$a_3a_4=ya_2~~\mbox{and}~~x_3x_4=zx_2.
$
All these together reduces \eqref{eq317317} to
\begin{align} &\frac{a_1D_{0}(\bar{a},\bar{p})}{qa_2(1-a_1y/q)}
\sum_{k=0}^{n-1}
q^{k}
\frac{(a_1^2/q,x_1^2;q)_k}{
(a_1y,x_1z;q)_{k}} \label{eq317} \\
&=1- \frac{(a_1^2/q,x_1^2;q)_n}{
(a_1y/q,x_1z;q)_n} -\frac{x_1D_{0}(\bar{x},\bar{q})(1-a_1^2/q)}{x_2(1-a_1y/q)(1-x_1z)}\sum_{k=0}^{n-1}
q^k\frac{(a_1^2,x_1^2;q)_k}{(a_1y,x_1zq;q)_{k}}.
\nonumber
\end{align}
Upon substituting  $D_{0}(\bar{a},\bar{p})$ and $D_{0}(\bar{x},\bar{q})$ in \eqref{eq317} and  noting that $a_3a_4=ya_2$, we get
\begin{align*}
&\frac{a_1(a_1-y)}{q-a_1y}\sum_{k=0}^{n-1}
q^{k}
\frac{(a_1^2/q,x_1^2;q)_k}{
(a_1y,x_1z;q)_{k}}\\
&\quad=1- \frac{(a_1^2/q,x_1^2;q)_n}{
(a_1y/q,x_1z;q)_n}-\frac{x_1(x_1-z)(q-a_1^2)}{(q-a_1y)(1-x_1z)}\sum_{k=0}^{n-1}
q^k\frac{(a_1^2,x_1^2;q)_k}{(a_1y,x_1zq;q)_{k}}.
\nonumber
\end{align*}
Finally, expressed in terms of $\mathbf{\Pi}_n(a_1,x_1,y,z)$, it takes the form of \eqref{eq317317}.
As it turns out, \eqref{xrma-1-xrma} follows from \eqref{eq317317} after letting $x_1z=q$ and then taking $n\to \infty$.
\qed

We think that \eqref{xrma-1-xrma} is novel to the literature,  after compared with  \cite[17.6.E13]{dlmf}. Far beyond this, on appealing to the contiguous relation \eqref{eq317317}, we may establish
\begin{tl}
\begin{align}
\sum_{k=0}^{\infty}
q^{k}
\frac{(a_1^2/q,x_1^2;q)_k}{
(a_1y,x_1z;q)_{k}}&=\bigg(\frac{x_1^2q}{a_1y}\bigg)^m
\frac{(z/x_1,a_1^2/q;q)_m}{(a_1/y,x_1z;q)_m}
\sum_{k=0}^{\infty}
q^{k}
\frac{(a_1^2q^{m-1},x_1^2;q)_k}{
(a_1y,x_1zq^{m};q)_{k}}\nonumber\\
&+\frac{q-a_1y}{a_1(a_1-y)}\sum_{k=0}^{m-1}
\bigg(\frac{x_1^2q}{a_1y}\bigg)^k
\frac{(z/x_1,a_1^2/q;q)_k}{(a_1q/y,x_1z;q)_k}\label{eq317317-new}\\
&-\frac{q}{a_1(a_1-y)}\frac{(x_1^2,a_1^2/q;q)_\infty}{(a_1y,x_1z;q)_\infty}
\sum_{k=0}^{m-1}
\bigg(\frac{x_1^2q}{a_1y}\bigg)^k\frac{(z/x_1;q)_k}{(a_1q/y;q)_k}.\nonumber
\end{align}
\end{tl}
\pf Actually, \eqref{eq317317} can be restated as the operator form
\begin{align}
\mathbf{\Pi}_n(a_1,x_1,y,z)&=C(a_1,x_1,y,z)
\sigma\big(\mathbf{\Pi}_n(a_1,x_1,y,z)\big)+D_n(a_1,x_1,y,z),\label{xrma-1-11-111}
\end{align}
where the operator
\[\sigma(a_1,x_1,y,z):=(a_1q^{1/2},x_1,y/q^{1/2},zq)\]
and the coefficients
\begin{align*}C(a_1,x_1,y,z)&:=
\frac{x_1(x_1-z)(a_1^2-q)}{a_1(a_1-y)(1-x_1z)},\\
D_n(a_1,x_1,y,z)&:=\frac{q-a_1y}{a_1(a_1-y)}
\left\{1-\frac{(a_1^2/q,x_1^2;q)_n}{(a_1y/q,x_1z;q)_n}\right\}.\end{align*}
Hence, \eqref{eq317317-new} comes out by applying Lemma \ref{xxx-xxx} to \eqref{xrma-1-11-111}
 and then taking the limit as $n\to \infty$.
\qed

\begin{tl}
For  $|x_1^2q/a_1y|<1$, it holds
\begin{align}
\sum_{k=0}^{\infty}
q^{k}
\frac{(a_1^2/q,x_1^2;q)_k}{
(a_1y,x_1z;q)_{k}}&=\frac{q-a_1y}{a_1(a_1-y)}\sum_{k=0}^{\infty}
\bigg(\frac{x_1^2q}{a_1y}\bigg)^k\frac{(z/x_1,a_1^2/q;q)_k}
{(a_1q/y,x_1z;q)_k}\label{eq317317-new-new}\\
&-\frac{q}{a_1(a_1-y)}\frac{(x_1^2,a_1^2/q;q)_\infty}{(a_1y,x_1z;q)_\infty}\sum_{k=0}^{\infty}
\bigg(\frac{x_1^2q}{a_1y}\bigg)^k\frac{(z/x_1;q)_k}{(a_1q/y;q)_k} .\nonumber
\end{align}
Further, for $x_1z=q$, it holds
\begin{align}
{}_2 \phi_1\left[\begin{array}{c}
a_1^2 / q, x_1^2 \\
a_1 y
\end{array} ; q, q\right]&+\frac{q}{a_1(a_1-y)}\frac{\big(a_1^2 / q,x_1^2; q\big)_{\infty}}
{(a_1 y,q; q)_{\infty}}{}_2 \phi_1\left[\begin{array}{c}
q, q/x_1^2 \\
a_1q/y
\end{array}; q, \frac{x_1^2q}{a_1y}\right]\nonumber\\
&=\frac{(q/a_1y,a_1x_1^2/y;q)_\infty}{(a_1/y,x_1^2q/a_1y;q)_\infty}.
\label{xrma-1-xrma5}
\end{align}
\end{tl}
\pf  It is clear that \eqref{eq317317-new-new} is the limiting case of \eqref{eq317317-new} as $m\to\infty$. Note that  $\lim_{m\to\infty}\big(x_1^2q/a_1y\big)^m=0$ for $|x_1^2q/a_1y|<1$.
As a result, for $x_1z=q$, we can reformulate \eqref{eq317317-new-new} as
\begin{align}
{}_2 \phi_1\left[\begin{array}{c}
a_1^2 / q, x_1^2 \\
a_1 y
\end{array} ; q, q\right]&+\frac{q-a_1y}{a_1(y-a_1)}{}_2 \phi_1\left[\begin{array}{c}
a_1^2/q, q/x_1^2\\
a_1q/y
\end{array} ; q, \frac{x_1^2q}{a_1y}\right]\nonumber\\
&=\frac{q}{a_1(y-a_1)}\frac{\big(a_1^2 / q,x_1^2; q\big)_{\infty}}
{(a_1 y,q; q)_{\infty}}{}_2 \phi_1\left[\begin{array}{c}
q, q/x_1^2 \\
a_1q/y
\end{array}; q, \frac{x_1^2q}{a_1y}\right].\label{xrma-1-xrma3}
\end{align}
In this case, by the $q$-Gauss summation formula \cite[(II. 8)]{10}, we find
$$
{}_2 \phi_1\left[\begin{array}{c}
a_1^2/q, q/x_1^2\\
a_1q/y
\end{array} ; q,
\frac{x_1^2q}{a_1y}\right]=
\frac{(q^2/a_1y,a_1x_1^2/y;q)_\infty}{(a_1q/y,x_1^2q/a_1y;q)_\infty}.
$$
Substituting this into \eqref{xrma-1-xrma3}, we have
\begin{align*}
{}_2 \phi_1\left[\begin{array}{c}
a_1^2 / q, x_1^2 \\
a_1 y
\end{array} ; q, q\right]&-
\frac{(q/a_1y,a_1x_1^2/y;q)_\infty}{(a_1/y,x_1^2q/a_1y;q)_\infty}\nonumber\\
&=\frac{q}{a_1(y-a_1)}\frac{\big(a_1^2 / q,x_1^2; q\big)_{\infty}}
{(a_1 y,q; q)_{\infty}}{}_2 \phi_1\left[\begin{array}{c}
q, q/x_1^2 \\
a_1q/y
\end{array}; q, \frac{x_1^2q}{a_1y}\right].
\end{align*}
This amounts to \eqref{xrma-1-xrma5}.
\qed
\section{Concluding remarks}
A quick glance over all that we have done in the foregoing sections together with  \cite{xuxrma}  inspires us to  wonder what  will Theorem \ref{firsteq} subject to two conditions
\begin{align*}
\left\{\begin{matrix}~a_1^2/p_1=x_1x_2Y\\
~a_2^2/p_2=x_1x_2/Y
\end{matrix}
\right. \quad\mbox{and}\quad \left\{\begin{matrix}~p_1=L_0/q_3\\
~p_2=L_0/q_4
\end{matrix}
\right.\qquad\mbox{become?}
\end{align*}
  As a preliminary result, we have now shown that
\begin{dl}\label{firsteq-simplified}With the same notation as Theorem \ref{firsteq}. Let $L=(p_1p_2p_3p_4)^{1/2}$ and further assume that $p_1p_2=q_1q_2$ and $a_1^2a_2^2=x_1^2x_2^2.$ Then, for any integer $n\geq 0$, it holds
\begin{align}
\sum_{k=0}^{n-1}g(k)\frac{(x_1^2;q_1)_k(x_2^2;q_2)_k}{(a_1a_2X;L/p_3)_k(a_1a_2/X;L/p_4)_k
}=\frac{(x_1^2;q_1)_n(x_2^2;q_2)_n}{
(a_1a_2X;L/p_3)_{n}(a_1a_2/X;L/p_4)_{n}
}-1,
\end{align}
where
\begin{align}
g(k):=(L/p_3)^k\frac{(a_1a_2X-x_1^2(q_1p_3/L)^k)(a_1a_2X-x_2^2(q_2p_3/L)^k)}{a_1a_2
(1-a_1a_2X(L/p_3)^{k})(X-a_1a_2(L/p_4)^k)}.
\end{align}
\end{dl}
In a broad range, the same  question about  \eqref{pppppp}  of Lemma \ref{type-i-i} is also of interest. Truly,  the answer to this problem  would become more complicated to figure out. By contrast,  any  reduction of \eqref{oooooo-1} seems highly possibility for good $q$-contiguous relations of the truncated series involved.
Further discussion on this topic will be made in our forthcoming paper.

\section*{Acknowledgement}
This  work was supported by the National Natural Science Foundation of
China [Grant No. 11971341].

\end{document}